\documentclass[12pt,a4paper]{article}
\usepackage[T2A]{fontenc}
\usepackage[cp1251]{inputenc}
\usepackage[russian]{babel}
\usepackage{amsmath,amssymb,amsthm,amsfonts,amscd}

\begin{document}
\sloppy
\date{}

\title{Geometric approach towards stable homotopy groups of
spheres. The Kervaire invariant}

\author{P.M.Akhmet'ev \thanks{This work was supported in part by the London Royal Society
(1998-2000), RFBR 08-01-00663,  INTAS 05-1000008-7805.}}


\newtheorem{theorem}{Теорема}[section]
\newtheorem*{main*}{Основная Теорема}
\newtheorem*{theorem*}{Теорема}
\newtheorem{lemma}[theorem]{Лемма}
\newtheorem{proposition}[theorem]{Предложение}
\newtheorem{corollary}[theorem]{Следствие}
\newtheorem{conjecture}[theorem]{Гипотеза}
\newtheorem{problem}[theorem]{Проблема}

\theoremstyle{definition}
\newtheorem{definition}[theorem]{Определение}
\newtheorem{remark}[theorem]{Замечание}
\newtheorem*{remark*}{Замечание}
\newtheorem*{example*}{Пример}
\def\Z{{\Bbb Z}}
\def\R{{\Bbb R}}
\def\RP{{\Bbb R}\!{\rm P}}
\def\N{{\Bbb N}}
\def\C{{\bf C}}
\def\A{{\bf A}}
\def\D{{\bf D}}
\def\O{{\bf O}}
\def\I{{\bf I}}

\def\fr{{\operatorname{fr}}}
\def\st{{\operatorname{st}}}
\def\mod{{\operatorname{mod}\,}}
\def\cyl{{\operatorname{cyl}}}
\def\dist{{\operatorname{dist}}}
\def\sf{{\operatorname{sf}}}
\def\dim{{\operatorname{dim}}}

\maketitle

\begin{abstract}
The notion of the geometrical $\Z/2 \oplus \Z/2$--control of
self-intersection of a skew-framed immersion and the notion of the
$\Z/2 \oplus \Z/4$-structure (the cyclic structure) on the
self-intersection manifold of a $\D_4$-framed immersion are
introduced. It is shown that a skew-framed immersion
$f:M^{\frac{3n+q}{4}} \looparrowright \R^n$, $0 < q <<n$ (in the
$\frac{3n}{4}+\varepsilon$-range) admits a geometrical $\Z/2
\oplus \Z/2$--control if the characteristic class of the
skew-framing of this immersion admits a retraction of the order
$q$, i.e. there exists a mapping $\kappa_0: M^{\frac{3n+q}{4}} \to
\RP^{\frac{3(n-q)}{4}}$, such that this composition $I \circ
\kappa_0: M^{\frac{3n+q}{4}} \to \RP^{\frac{3(n-q)}{4}} \to
\RP^{\infty}$ is the characteristic class of the skew-framing of
$f$. Using the notion of $\Z/2 \oplus \Z/2$-control  we prove that
for a sufficiently great $n$, $n=2^l-2$, an arbitrary immersed
$\D_4$-framed manifold  admits in the regular cobordism class
(modulo odd torsion) an immersion with a $\Z/2 \oplus
\Z/4$-structure.
In the last section we present an approach toward  the Kervaire Invariant One Problem.
\end{abstract}

\section{Self-intersection of immersions and Kervaire Invariant}
The Kervaire Invariant One Problem is an open problem in Algebraic
topology, for algebraic approach see [B-J-M], [C-J-M]. We will
consider a geometrical approach;
this approach is based on results by P.J.Eccles, see [E1]. For a
geometrical approach see also [C1],[C2].

Let $f: M^{n-1} \looparrowright \R^n$, $n= 2^l -2$, $l>1$, be a
smooth (generic) immersion of codimension 1. Let us denote by $g:
N^{n-2} \looparrowright \R^n$ the immersion of self-intersection
manifold.

\subsubsection*{Definition 1}
The Kervaire invariant of $f$ is defined as
$$ \Theta(f) = <w_2^{\frac{n-2}{2}}; [N^{n-2}] >, $$
where $w_2 = w_2(N^{n-2})$ is the normal Stiefel-Whitney of
$N^{n-2}$.

\[  \]
The Kervaire invariant is an invariant of the regular cobordism
class of the immersion $f$. Moreover, the Kervaire invariant is
a well-defined homomorphism
$$ \Theta: Imm^{sf}(n-1,1) \to  \Z/2. \eqno(1) $$

The normal bundle $\nu(g)$ of the immersion $g: N^{n-2}
\looparrowright \R^n$ is a 2-dimensional bundle over $N^{n-2}$
equipped with a $\D_4$--framing. The classifying mapping $\eta:
N^{n-2} \to K(\D_4,1)$ of this bundle is well-defined. The
$\D_4$-structure of the normal bundle or the $\D_4$--framing is
the prescribed reduction of the structure group of the normal
bundle of the immersion $g$ to the group $\D_4$ corresponding to
the mapping $\eta$. The pair $(g,\eta)$ represents an element in
the cobordism group $Imm^{\D_4}(n-2,2)$. The homomorphism
$$ \delta: Imm^{sf}(n-1,1) \to Imm^{\D_4}(n-2,2) \eqno(2)$$
is well-defined.

Let us recall that the cobordism group $Imm^{sf}(n-k,k)$
generalizes the group $Imm^{sf}(n-1,1)$. This group is defined as
the cobordism group of triples $(f,\Xi,\kappa)$, where $f: M^{n-k}
\looparrowright \R^n$ is an immersion with the prescribed
isomorphism $\Xi: \nu(g) \cong k \kappa$, called a skew-framing,
$\nu(f)$ is the normal bundle of $f$, $\kappa$ is the given line
bundle over $M^{m-k}$ with the characteristic class $w_1(\kappa)
\in H^1(M^{m-k};\Z/2)$. The cobordism relation of triples is
standard.

The generalization of the group $Imm^{\D_4}(n-2,2)$ is following.
Let us define the cobordism groups $Imm^{\D_4}(n-2k,2k)$. This
group $Imm^{\D_4}(n-2k,2k)$ is represented by triples
$(g,\Xi,\eta)$, where $g: N^{n-2k} \looparrowright \R^n$ is an
immersion, $\Xi$ is a dihedral $k$-framing, i.e. the prescribed
isomorphism $\Xi: \nu_g \cong k \eta$, where $\eta$ is a 2-dimensional
bundle over $N^{n-2k}$. The characteristic mapping of the bundle
$\eta$ is denoted also by $\eta: N^{n-2k} \to K(\D_4,1)$. The
mapping $\eta$ is the characteristic mapping for the bundle
$\nu_g$, because $\nu_g \cong k \eta$.

Obviously, the Kervaire homomorphism (1) is defined as the
composition of the homomorphism (2) with a homomorphism
$$ \Theta_{\D_4} : Imm^{\D_4}(n-2,2) \to \Z/2. \eqno(3) $$
The homomorphism (3) is called the Kervaire invariant for
$\D_4$-framed immersed manifolds.

The Kervaire homomorphisms are defined in a more general situation
by a straightforward generalization of the homomorphisms (1) and
(3):
$$ \Theta^k: Imm^{sf}(n-k,k) \to \Z/2, \eqno(4a) $$
$$ \Theta^k_{\D_4} : Imm^{\D_4}(n-2k,2k) \to \Z/2, \eqno(4b) $$
(for $k=1$ the new homomorphism coincides with the homomorphism
(3) defined above) and the following diagram
$$
\begin{array}{ccccc}
Imm^{sf}(n-1,1) & \stackrel {\delta}{\longrightarrow} &
Imm^{\D_4}(n-2,2) & \stackrel{\Theta_{\D_4}}{\longrightarrow} & \Z/2  \\
\downarrow J^k & & \downarrow J^k_{\D_4}  &  &  \vert \vert \\
Imm^{sf}(n-k,k) & \stackrel{\delta^k}{\longrightarrow} & Imm^{\D_4}(n-2k,2k) &
\stackrel{\Theta_{\D_4}^k}{\longrightarrow} & \Z/2  \\
\end{array} \eqno(5)
$$
is commutative. The homomorphism $J^k$ ($J^k_{\D_4}$) is determined by the regular cobordism class
of the restriction of the given immersion $f$ ($g$) to the submanifold in $M^{n-1}$ ($N^{n-2}$) dual to
$w_1(\kappa)^{k-1} \in H^{k-1}(M^{n-1};\Z/2)$ ($w_2(\eta)^{k-1} \in H^{2k-2}(N^{n-2};\Z/2)$).

Let $(g,\Xi,\eta)$ be a $\D_4$-framed (generic) immersion in the
codimension $2k$. Let $h: L^{n-4k} \looparrowright \R^n$ be the
immersion of the self-intersection (double points) manifold of
$g$. The normal bundle $\nu_h$ of the immersion $h$ is decomposed
into a direct sum of $k$ isomorphic copies of a 4-dimensional bundle
$\zeta$ with the structure group $\Z/2 \int \D_4$. This
decomposition is given by the isomorphism $\Psi: \nu_h \cong k \zeta$.
The bundle $\nu_h$ itself is classified by the mapping $\zeta:
L^{n-4k} \to K(\Z/2 \int \D_4,1)$.

All the triples $(h,\zeta,\Psi)$ described above (we do not assume
that a triple is realized as the double point manifold for a
$\D_4$-framed immersion) up to the standard cobordism relation
form the cobordism group $Imm^{\Z/2 \int \D_4}(n-4k,4k)$. The
self-intersection of an arbitrary $\D_4$-framed immersion is a
$\Z/2 \int \D_4$-framed immersed manifold and the cobordism class
of this manifold well-defines the natural homomorphism
$$ \delta_{\D_4}^k : Imm^{\D_4}(n-2k,2k) \to Imm^{\Z/2 \int \D_4}(n-4k,4k). \eqno(6) $$

The subgroup $\D_4 \oplus \D_4 \subset \Z/2 \int \D_4$ of 
index 2 induces the double cover $\bar L^{n-4k} \to L^{n-4k}$.
This double cover corresponds with the canonical double cover over the
double point manifold.

Let $\bar \zeta: \bar L^{n-4k} \to K(\D_4,1)$ be the classifying
mapping induced by the projection homomorphism $\D_4 \oplus \D_4
\to \D_4$ to the first factor. Let $\bar \zeta \to L^{n-4k}$ be
the 2-dimensional $\D_4$--bundle defined as the pull-back of the
universal 2-dimensional bundle with respect to the classifying
mapping $\bar \zeta$.

\subsubsection*{Definition 2}
The Kervaire invariant $\Theta_{\Z/2 \int \D_4}^k: Imm^{\Z/2 \int \D_4}(n-4k,4k) \to \Z/2$ for a $\Z/2
\int \D_4$-framed immersion $(h,\Psi,\zeta)$ is defined by the
following formula:
$$ \Theta_{\Z/2 \int \D_4}^k(h,\Psi,\zeta) = <w_2(\bar \eta)^{\frac{n-4k}{2}};[L^{n-4k}]>. $$
\[  \]

This new invariant is a homomorphism $\Theta_{\Z/2 \int \D_4}^k:
Imm^{\Z/2 \int \D_4}(n,n-4k) \to \Z/2$ included into the following
commutative diagram:

$$
\begin{array}{ccc}

Imm^{\D_4}(n-2k,2k) & \stackrel{\Theta_{\D_4}}{\longrightarrow} & \Z/2  \\
\downarrow  \delta_{\D_4}^k  &  &  \vert \vert \\
Imm^{\Z/2 \int \D_4}(n-4k,4k) &
\stackrel{\Theta_{\Z/2 \int \D_4}^k}{\longrightarrow} & \Z/2.  \\
\end{array} \eqno(7)
$$

Let us formulate the first main results of the paper. In 
section 2 the notion of $\Z/2 \oplus \Z/2$-control
($\I_b$--control) on self-intersection of a skew-framed immersion
is considered. Theorem 1 (for the proof see section 3) shows
that under a natural restriction of dimensions the property of
$\I_b$-control holds for an immersion in the regular cobordism
class modulo odd torsion.

In section 4 we formulate a notion of $\Z/2 \oplus
\Z/4$--structure (or an $\I_4$--structure, or a cyclic structure)
of a $\D_4$-framed immersion. In section 5 we prove Theorem 2. We prove under a
natural restriction of dimension that an arbitrary $\D_4$-framed
$\I_b$-controlled immersion admits in the regular homotopy class
an immersion with a cyclic structure. For such an immersion 
Kervaire invariant is expressed in terms of $\Z/2 \oplus
\Z/4$--characteristic numbers of the self-intersection manifold.
The proof (based on the two theorems from [A2] (in Russian))
of the Kewrvaire Invariant One Problem is in section 6.

The author is grateful to Prof. M.Mahowald (2005) and Prof.
R.Cohen (2007) for discussions, to Prof. Peter Landweber for the help with the English translation, 
and to  Prof. A.A.Voronov for the
invitation to Minnesota University in (2005). 

This paper was started in 1998 at the Postnikov Seminar. This
paper is dedicated to the memory of Prof. Yu.P.Soloviev.

\section{Geometric Control of self-intersection manifolds of
skew-framed immersions} In this and the remining sections of the paper  by
$Imm^{sf}(n-k,k)$, $Imm^{\D_4}(n-2k,2k)$, $Imm^{\Z/2 \int
\D_4}(n-4k,4k)$, etc., we will denote not the cobordism groups
themselves, but the 2-components of these groups. In case the
first argument (the dimension of the immersed manifold) is
strictly positive, all the groups are finite 2-group.

Let us recall that the dihedral group $\D_4$ is given by the
representation (in terms of generators and relations) $\{a,b\vert a^4 = b^2 = e, [a,b]=a^2\}$. This
group is a subgroup of the group $O(2)$ of isometries of the
plane with the base $\{f_1,f_2\}$ that keeps the pair of lines
generated by the vectors of the base. The element $a$ corresponds
to the rotation of the plane through the angle $\frac{\pi}{2}$. The
element $b$ corresponds to the reflection of the plane with the
axis given by the vector $f_1 + f_2$.

Let $\I_b(\Z/2 \oplus \Z/2)= \I_b \subset \D_4$ be the subgroup
generated by the elements $\{a^2,b\}$. This is an elementary
$2$-group of rank 2 with two generators. These are the transformations of the plane
that preserve each line $l_1$, $l_2$ generated by the vectors
$f_1+f_2$, $f_1-f_2$ correspondingly. The cohomology group
$H^1(K(\I_b,1);\Z/2)$ is the elementary 2-group with two
generators. The first (second) generator of this group detects the
reflection of the line $l_2$  (of the line $l_1$) correspondingly.
 The generators of the cohomology
group will be denoted by $\tau_1$, $\tau_2$ correspondingly.

\subsubsection*{Definition 3}

We shall say that a skew-framed immersion $(f,\Xi)$, $f: M^{n-k}
\looparrowright \R^n$ has self-intersection of type $\I_b$, if
the double-points manifold $N^{n-2k}$ of $f$ is a $\D_4$-framed
manifold that admits a reduction of the structure group $\D_4$ of
the normal bundle to the subgroup $\I_b \subset \D_4$.
\[  \]

Let us formulate the following conjecture.

\subsubsection*{Conjecture}

For an arbitrary $q > 0$, $q=2 (mod 4)$, there exists a positive
integer $l_0=l_0(q)$, such that for an arbitrary $n = 2^l-2$,
$l>l_0$ an arbitrary element $a \in
Imm^{sf}(\frac{3n+q}{4},\frac{n-q}{4})$ is stably regular
cobordant to a stably skew-framed immersion with $\I_b$-type of
self-intersection (for the definition of stable framing see [E2], of
stable skew-framing see [A1]).
\[  \]

Let us formulate and prove a weaker result toward the Conjecture.
We start with the following definition.

Let $\omega: \Z/2 \int \D_4 \to \Z/2$ be the epimorphism defined
as the composition $\Z/2 \int \D_4 \subset \Z/2 \int \Sigma_4 \to
\Sigma_4 \to \Z/2$, where $\Sigma_4 \to \Z/2$ is the parity of a
permutation. Let $\omega^{!}: Imm^{\Z/2 \int \D_4}(n-4k,4k) \to
Imm^{Ker \omega}(n-4k,4k)$ be the transfer homomorphism with
respect to the kernel of the epimorphism $\omega$.

Let $P$ be a polyhedron with $dim(P) < 2k-1$, $Q \subset P$ be a
subpolyhedron with $dim(Q)=dim(P)-1$, and let $P \subset \R^n$ be an
embedding. Let us denote by  $U_P$ the regular neighborhood of $P
\subset \R^n$ of the radius $r_P$ and by $U'_Q$ the regular
neighborhood of $Q \subset \R^n$ of the radius $r_Q$, $r_Q > r_P$.
Let us denote $U_Q = U_P \cap U'_Q$.

The boundary $\partial U_P$ of the neighborhood $U_P$ is a
codimension one submanifold in $\R^n$. This manifold $\partial
U_P$ is a union of the two manifolds with boundaries
 $V_Q \cup_{\partial} V_P$,  $V_Q = U_Q \cap
\partial U_P$, $V_P = \partial U_P \setminus U_Q$
along the common boundary $\partial V_Q = \partial V_P$.

Let us assume that the two cohomology classes $\tau_{Q,1} \in
H^1(Q;\Z/2)$, $\tau_{Q,2} \in H^1(Q;\Z/2)$ are given. The
projection $U_Q \to Q$ of the neighborhood on the central
submanifold determines the cohomology classes $\tau_{U_Q,1},
\tau_{U_Q,2} \in H^1(U_Q;\Z/2)$ as the inverse images of the
classes $\tau_{Q,1}, \tau_{Q,2}$ correspondingly.

Let $(g,\Xi_N,\eta)$, $dim(N)=n-2k$ be a  $\D_4$--framed generic immersion, $n-4k > 0$, and
 $g(N^{n-2k}) \cap \partial U_P$ be an immersed submanifold in $U_Q \subset \partial U_P$.
 Let us denote
$g(N^{n-2k}) \setminus (g(N^{n-2k}) \cap (U_P))$ by
$N^{n-2k}_{int}$, and the complement $N^{n-2k} \setminus
N^{n-2k}_{int}$ by $N^{n-2k}_{ext}$. The manifolds
$N^{n-2k}_{ext}$, $N^{n-2k}_{int}$ are submanifolds in $N^{n-2k}$
of codimension 0 with the common boundary, this boundary is
denoted by $N_Q^{n-2k-1}$. The self-intersection manifold of $g$
is denoted by $L^{n-4k}$. By the dimensional reason ($n-4k=q<<n$)
$L^{n-4k}$ is a submanifold in $\R^n$, parameterized by an
embedding $h$,
 equipped by the $\Z/2 \int \D_4$-framing of the normal bundle denoted by $(\Psi, \zeta)$.
The triple $(h,\Psi,\zeta)$ determines an element in the cobordism
group $Imm^{\Z/2 \int \D_4}(n-4k,4k)$.

\subsubsection*{Definition 4}
We say that the $\D_4$--framed
immersion $g$ is an $\I_b$--controlled immersion if the following
conditions hold:

--1. The structure group of the $\D_4$--framing $\Xi_N$ restricted
to the submanifold (with boundary) $g(N^{n-2k}_{ext})$ is
reduced to the subgroup $\I_b \subset \D_4$ and  the cohomology
classes $\tau_{U_Q,1}, \tau_{U_Q,2} \in H^1(U_Q;\Z/2)$ are mapped
to the generators $\tau_1, \tau_2 \in H^1(N_Q^{n-2k-1};\Z/2)$ of
the cohomology of the structure group of this $\I_b$-framing by
the immersion $g \vert_{N^{n-2k-1}_Q} : N^{n-2k-1}_Q
\looparrowright
\partial(U_Q) \subset U_Q$.

--2. The restriction of the immersion $g$ to the submanifold
 $N_Q^{n-2k-1} \subset N^{n-2k}$ is an embedding
$g \vert_{N_Q^{n-2k-1}} : N_Q^{n-2k-1} \subset \partial U_Q$, and
the decomposition $L^{n-4k} =L^{n-4k}_{int} \cup L^{n-4k}_{ext}
\subset (U_P \cup \R^n \setminus U_P)$ of the self-intersection
manifold of $g$ into two (probably, non-connected) $\Z/2 \int
\D_4$-framed components is well-defined. The manifold
$L^{n-4k}_{int}$ is a submanifold in $U_P$ and the triple
$(L^{n-4k}_{int},\Psi_{int},\zeta_{int})$ represents an element
in $Imm^{Ker \omega}(n-4k,4k)$ in the image of the homomorphism
$\omega^!: Imm^{\Z/2 \int \D_4}(n-4k,4k) \to Imm^{Ker
\omega}(n-4k,4k)$.


\[  \]

\subsubsection*{Definition 5}
Let $(f,\Xi_M,\kappa) \in Imm^{sf}(n-k,k)$ be an arbitrary
element, where
 $f: M^{n-k} \looparrowright \R^n$ is an immersion of codimension $k$ with the characteristic class
$\kappa \in H^1(M^{n-k};\Z/2)$ of the skew-framing $\Xi_M$. We say that
the pair $(M^{n-k},\kappa)$ admits a retraction of order $q$, if the
mapping $\kappa : M^{n-k} \to \RP^{\infty}$ is represented by the
composition $\kappa = I \circ \bar \kappa : M^{n-k} \to
\RP^{n-k-q-1} \subset \RP^{\infty}$. The element $[(f, \Xi_M,
\kappa)]$ admits a retraction of order $q$, if in the
cobordism class of this skew-framed immersion there exists a
triple $(M'^{n-k}, \Xi_{M'}, \kappa')$ that admits a retraction of 
order $q$.

\[  \]

\subsubsection*{Theorem 1}

Let $q = q(l)$ be a positive integer,
 $q=2(mod 4)$. Let us assume that an element $\alpha \in Imm^{sf}(\frac{3n+q}{4},\frac{n-q}{4})$
 admits a retraction of the order $q$ and $3n-12k-4>0$.  Then the element $\delta(\alpha) \in
 Imm^{\D_4}(n-2k,2k)$, $k=\frac{n-q}{4}$,
 is represented by a $\D_4$-framed immersion
  $[(g,\Psi_N,\eta)]$ with  $\I_b$-control.
\[  \]

\section{Proof of  Theorem 1}

Let us denote $n-k-q-1 = 3k-1$ by $s$. Let
 $d: \RP^{s} \to \R^n$ be a generic mapping. We denote the self-intersection points of $d$
 (in the target space) by $\Delta(d)$ and the singular points of $d$ by $\Sigma(d)$.

Let us recall a classification of singular points of generic
mappings $\RP^{s} \to \R^n$ in the case $4s < 3n$, for details see
[Sz]. In this range generic mappings have no quadruple points. The
singular values (in the target space) are of the following two
types:

-- a closed manifold $\Sigma^{1,1,0}$;

-- a singular manifold $\Sigma^{1,0}$ (with singularities of the
type $\Sigma^{1,1,0}$).

The multiple points are of the multiplicities 2 and 3. The set of triple
points form a manifold with boundary and with corners on the
boundary. These "corner" \- singular points on the boundary of the
triple points manifold coincide  with the manifold
$\Sigma^{1,1,0}$. The regular part of boundary of triple points is
a submanifold in $\Sigma^{1,0}$.

The double self-intersection points form a singular submanifold in $\R^n$ with the
boundary $\Sigma^{1,0}$. This submanifold is not generic. After an
arbitrary small alteration the double points manifold becomes a
submanifold in $\R^n$ with boundary and with corners on the
boundary of the type $\Sigma^{1,1,0}$.

Let  $U_{\Sigma}$ be a small regular neighborhood of the radius
$\varepsilon_1$ of the singular submanifold $\Sigma^{1,0}$. Let
$U_{\Delta}$ be a small regular neighborhood of the same radius of
the submanifold $\Delta(d)$ (this submanifold is immersed with
singularities on the boundary). The inclusion $U_{\Sigma} \subset
U_{\Delta}$ is well-defined.

Let us consider a regular submanifold in $\Delta$ obtained by 
excising a small regular neighborhood of the boundary.  This
immersed manifold with boundary will be denoted by $\Delta^{reg}$.
The (immersed) boundary $\partial \Delta^{reg}$ will be denoted by
$\Sigma^{reg}$. We will consider the pair of  regular
neighborhoods $U^{reg}_{\Sigma} \subset U^{reg}_{\Delta}$ of the
pair $\Sigma^{reg} \subset \Delta^{reg}$ of the radius
$\varepsilon_2$, $\varepsilon_2 << \varepsilon_1$. Because $2
dim(\Delta^{reg})<n$, after a small perturbation the manifold
$\Delta^{reg}$ is a submanifold  in $U^{reg}_{\Delta}$.

Let $(f_0,\Xi_0,\kappa)$, $f_0: M^{n-k} \looparrowright \R^n$,
$n-k=\frac{3n + q}{4}$ be a skew-framed immersion in the
cobordism class $\alpha$. We will construct an immersion $f:
M^{n-k} \looparrowright \R^n$ in the regular homotopy class of
$f_0$ by the following construction.

Let  $\kappa_0: M^{n-k} \to \RP^{s}$ be a retraction of order
$q$.  Let $f: M \looparrowright \R^n$ be an immersion in the
regular homotopy class of $f_0$ under the condition $dist(d \circ
\kappa_0, f_0) < \varepsilon_3$. The caliber $\varepsilon_3$ of
the approximation is given by the following inequality:
$\varepsilon_3 << \varepsilon_2$.

Let $g_1: N^{n-2k} \looparrowright \R^n$ be the immersion,
parameterizing the double points of $f$. The immersion $g_1$ is not
generic. After a small perturbation of the immersion $g_1$ with
the caliber $\varepsilon_3$  we obtain a generic immersion $g_2:
N^{n-2k} \looparrowright \R^n$.

The immersed submanifold $g_2(N^{n-2k})$ is divided into two
submanifolds $g_2(N^{n-2k}_{int})$, $g_2(N^{n-2k}_{ext})$ with the
common boundary $g_2(\partial N^{n-2k}_{int}) = g_2(\partial
N^{n-2k}_{ext})$ denoted by $g_2(N_Q^{n-2k-1})$. The manifold
$g_2(N^{n-2k}_{int})$ is defined as the intersection of the
immersed submanifold $g_2(N^{n-2k})$ with the neighborhood
$U^{reg}_{\Delta}$.  The manifold $g_2(N^{n-2k}_{ext})$ is defined
as the intersection of the immersed submanifold $g_2(N^{n-2k})$
with the complement  $\R^n \setminus (U^{reg}_{\Delta})$. We will
assume that $g_2$ is regular along $\partial U^{reg}_{\Delta}$.
Then $g_2(N^{n-2k}_Q)$ is an immersed submanifold in $\partial
U^{reg}_{\Delta}$. By construction the structure group $\D_4$
of the normal bundle of the immersed manifold
$g_2(N^{n-2k}_{ext})$ admits a reduction to the subgroup $\I_b
\subset \D_4$.

Let us denote by $L^{n-4k}$ the self-intersection manifold of the
immersion $g_2$. This manifold is embedded into $\R^n$ by $h:
L^{n-4k} \subset \R^n$. The normal bundle of this embedding $h$ is
equipped with a $\Z/2 \int \D_4$-framing denoted by $\Psi_L$ and the
characteristic class of this framing is denoted by $\zeta_L$. By
the analogous construction the manifold $L^{n-4k}$ is decomposed
as the union of the two manifolds over a common boundary, denoted
by $\Lambda$: $L^{n-4k}=L^{n-4k}_{ext} \cup_{\Lambda}
L^{n-4k}_{int}$. The manifold (with boundary) $L^{n-4k}_{int}$ is
embedded  by $h$ into $U^{reg}_{\Delta}$, the manifold
$L^{n-4k}_{ext}$ (with the same boundary) is embedded in the
complement $\R^n \setminus U^{reg}_{\Delta}$. The common boundary
$\Lambda$ is embedded into $\partial U^{reg}_{\Delta}$.

The manifold $L^{n-4k}$ is a $\Z/2 \int \D_4$-framed submanifold
in $\R^n$. Let us describe the reduction of the structure group of
this manifold to a corresponding subgroup in $\Z/2 \int \D_4$. We
will describe the subgroups
 $\I_{2,j}(\Z/2 \oplus \D_4)
\subset \Z/2 \int \D_4$, $j=x,y,z$.  We will
describe the transformations of $\R^4$ in the standard base
$(f_1,f_2,f_3,f_4)$ determined by generators of the groups.

Let us consider
the subgroup $\I_{2,x}$. The generator  $c_x$ (a generator will be
equipped with the index corresponding to the subgroup) defines the
transformation of the space by the following formula:
$c_x(f_1)=f_3$, $c_x(f_3)=f_1$, $c_x(f_2)=f_4$, $c_x(f_4)=f_2$.

For the generator $a_x$ (of the order 4) the transformation is
the following: $a_x(f_1)=f_2$, $a_x(f_2)=-f_1$, $a_x(f_3)=f_4$,
$a_x(f_4)=-f_3$. The generator $b_x$ (of order 2)  defines the
transformation of the space by the following formula:
 $b_x(f_1)=f_2$, $b_x(f_2)=f_1$, $b_x(f_3)=f_4$, $b_x(f_4)=f_3$.
From this formula the subgroup $\D_4 \subset \D_4 \oplus \Z/2$ is
represented by transformations that preserve the subspaces $(f_1,f_2)$,
$(f_3,f_4)$. The generator of the cyclic subgroup $\Z/2 \subset
\D_4 \oplus \Z/2$ permutes these planes.

The  subgroups $\I_{2,y}$ and $\I_{2,x}$ are conjugated by the
automorphism $OP: \Z/2 \int \D_4 \to \Z/2 \int \D_4$ given in the
standard base by the following formula: $f_1 \mapsto f_1$, $f_2 \mapsto
f_3$, $f_3 \mapsto f_2$, $f_4 \mapsto f_4$. Therefore the generator $c_y
\in \I_{2,y}$ is determined by the following transformation:
$c_y(f_1)=f_2$, $c_y(f_2)=f_1$, $c_y(f_3)=f_4$, $c_y(f_4)=f_3$.
The generator $a_y$ (of the order 4) is given by $a_y(f_1)=f_3$,
$a_y(f_3)=-f_1$, $a_y(f_2)=f_4$, $a_y(f_4)=-f_2$. The generator
$b_y$ (of the order 2) is given by $b_y(f_1)=f_3$, $b_y(f_3)=f_1$,
$b_y(f_2)=f_4$, $b_y(f_4)=f_2$.

Let us describe the subgroup $\I_{2,z}$. In this case the
generator $c_z$  defines the transformation of the space by the
following formula: $c_z(f_i)=-f_i$, $i=1,2,3,4$.

For the generator $a_z$ (of order 4) the transformation is the
following: $a_z(f_1)=f_2$, $a_z(f_2)=f_3$, $a_z(f_3)=f_4$,
$a_z(f_4)=f_1$. The generator $b_x$ (of the order 2)  defines the
transformation of the space by the following formula:
$b_z(f_1)=f_2$, $b_z(f_2)=f_1$, $b_z(f_3)=f_4$, $b_z(f_4)=f_3$.

Obviously, the restriction of the epimorphism $\omega: \Z/2 \int
\D_4 \to \Z/2$ to the subgroups $\I_{2,x}, \I_{2,y} \subset \Z/2
\int \D_4$ is trivial and the restriction of this homomorphism to
the subgroup $\I_{2,z}$ is non-trivial.

The subgroup $\I_3 \subset \I_{2,x}$ is defined as the subgroup
with the generators $c_x, b_x, a_x^2$. This is an index 2
subgroup isomorphic to the group $\Z/2^3$. The image of this
subgroup in $\Z/2 \int \D_4$ coincides with the intersection of
arbitrary pair of subgroups $\I_{2,x}$, $\I_{2,y}$, $\I_{2,z}$.
The subgroup $\I_3 \subset \I_{2,y}$ is generated by $c_y, b_y,
a^2_y$. Moreover, one has $c_y=b_x$, $b_y=c_x$, $a^2_y=a^2_x$. It
is easy to check that the following relations hold: $c_z=a^2_x$,
$a^2_z=c_x=b_y$, $b_z=b_x=c_y$. Therefore $Ker(\omega
\vert_{\I_2,z})$ coincides with the subgroup $\I_3 \subset
\I_{2,z}$.

The subgroups   $\I_{2,x}, \I_{2,y}, \I_{2,z}, \I_3$ in $\Z/2 \int
\D_4$ are well-defined. There is a natural projection $\pi_b: \I_3
\to \I_b$.

We will also consider the subgroup $\I_{2,x \downarrow} \subset
\Z/2 \int \D_4$ from geometrical considerations. This subgroup is a quadratic extension of the
subgroup $\I_{2,x}$ such that  $\I_{2,x}=Ker \omega \vert_{\I_{2,x
\downarrow}} \subset \I_{2,x\downarrow}$. An algebraic definition
of this group will not be required.

In the following lemma we will describe the structure group of the
framing of the triad $(L^{n-4k}_{int} \cup_{\Lambda}
L^{n-4k}_{ext})$. The framings of the spaces of the triad  will be denoted by $(\Psi_{\int}
\cup_{\Psi_\Lambda} \cup \Psi_{ext}, \zeta_{int}
\cup_{\zeta_\Lambda} \cup \zeta_{ext})$.

\subsubsection*{Lemma 1}

There exists a generic regular deformation $g_1 \to g_2$ of the
caliber $3 \varepsilon_3$ such that the immersed manifold
$g_2(N^{n-2k}_{ext})$ admits a reduction of the structure group
of the $\D_4$-framing to the subgroup $\I_b \subset \D_4$. The
manifold $L^{n-4k}_{int}$ is divided into the disjoint union of
the two manifolds (with boundaries) denoted by $(L^{n-4k}_{int,x
\downarrow}, \Lambda_{x \downarrow})$, $(L^{n-4k}_{int,y},
\Lambda_{y})$.

1. The structure group of the framing $(\Psi_{int,x \downarrow},
\Psi_{\Lambda_{x \downarrow}})$ for the submanifold (with
boundary) $(L^{n-4k}_{int,x \downarrow},\Lambda_{x \downarrow})$
is reduced to the subgroups $(\I_{2,x\downarrow},\I_{2,z})$.  (In
particular, the 2-sheeted cover over $L^{n-4k}_{int,x
\downarrow}$, classified by $\omega$ (denoted by $\tilde
L^{n-4k}_{int,x} \to L^{n-4k}_{int,x \downarrow}$)   is,
generally speaking, a non-trivial cover.)

2. The structure group of the framing $(\Psi_{int,y},
\Psi_{\Lambda})$  for the submanifold (with boundary)
$(L^{n-4k}_{int,y},\Lambda_y)$ is reduced to the subgroup
$(\I_{2,y},\I_3)$. (In particular, the 2-sheeted cover $\tilde
L^{n-4k}_{int,y} \to L^{n-4k}_{int,y}$ classified by $\omega$,
 is  the trivial cover.) Moreover, the
double covering $\tilde L^{n-4k}_{x}$ over the component
$L^{n-4k}_{x \downarrow}$ is naturally diffeomorphic to $\tilde
L^{n-4k}_y$ and this diffeomorphism  agrees with the restriction
of the automorphism $OP: \Z/2 \int \D_4 \to \Z/2 \int \D_4$ on the
subgroup $\I_{2,x}$, $OP(\I_{2,x})=\I_{2,y}$.

3.  The structure group of the framing $(\Psi_{ext}, \zeta_{ext})$
for the submanifold (with boundary) $h(L^{n-4k}_{ext},
\Lambda^{n-4k}) \subset (\R^n \setminus U^{reg}_{\Delta},
\partial(U^{reg}_{\Delta}))$ is reduced to the subgroup
 $\I_{2,z}$. (In particular, the
2-sheeted cover $\tilde L^{n-4k}_{ext} \to L^{n-4k}_{ext}$
classified by $\omega$, is, generally speaking, a nontrivial
cover.)

\[  \]

\subsubsection*{Proof of Lemma 1}

Components of the self-intersection manifold $g_1(N^{n-2k})
\setminus (g_1(N^{n-2k}) \cap U_{\Sigma})$ (this manifold is
formed by double points $x \in g_1(N^{n-2k}), x \notin U_{\Sigma}$
with inverse images $\bar x_1, \bar x_2 \in M^{n-k}$) are
classified by the following two types.

Type 1.
The points $\kappa(\bar x_1)$,
$\kappa(\bar x_2)$ in $\RP^s$ are $\varepsilon_2$-close.

Type 2. The distances between the points $\kappa(\bar x_1)$,
$\kappa(\bar x_2)$ in $\RP^s$ are  greater then the caliber
$\varepsilon_2$ of the regular approximation. Points of this type
belong to the regular neighborhood $U_{\Delta}$ (of the radius
$\varepsilon_1$).

Let us classify components of the triple self-intersection
manifold $\Delta_3(f)$ of the immersion $f$. The a priori classification of components is the
following.

A point $x \in
\Delta_3(f)$ has  inverse images $\bar x_1, \bar x_2, \bar x_3$ in
$M^{n-k}$.

Type 1. The images  $\kappa(\bar x_1), \kappa(\bar x_2),
\kappa(\bar x_3)$ are $\varepsilon_2$-close  in $\RP^s$.

Type 2. The images $\kappa(\bar x_1), \kappa(\bar x_2)$ are
$\varepsilon_2$-close  in $\RP^s$ and the distance between the
images $\kappa(\bar x_3)$ and $\kappa(\bar x_1)$ (or $\kappa(\bar
x_2)$) are greater than the caliber $\varepsilon_2$ of the
approximation.

Type 3. The pairwise distances between the points 
$\kappa(\bar x_1), \kappa(\bar x_2), \kappa(\bar x_3)$  greater than the
caliber $\varepsilon_2$ of the approximation.

By a general position argument the component of the type 3  does
not intersect  $d(\RP^s)$. Therefore the immersion $f$ can be
deformed by a small $\varepsilon_2$-small regular homotopy inside
the $\varepsilon_3$-regular neighborhood
 of the regular part of $d(\RP^s)$
such that after this regular homotopy $\Delta_3(f)$ is contained in
the complement of $U^{reg}_{\Delta}$. The codimension of the
submanifold $\bar \Delta_2(d) \subset \RP^s$ is equal to
$n-3k+1=q+k+1$ and  greater then $dim(\Delta_3(f)) = n-3k$. By
analogical arguments the component of triple points of the type 1
is outside $U^{reg}_{\Delta}$.

 Let us classify components of the quadruple self-intersection manifold
$\Delta_4(f)$ of the immersion $f$.  A point $x \in \Delta_4(f)$
has  inverse images $\bar x_1, \bar x_2, \bar x_3, \bar x_4$ in
$M^{n-k}$. The a priori classification is the following.

Type 1. The images $\kappa(\bar x_1), \kappa(\bar x_2)$ are
$\varepsilon_2$-close in $\RP^s$ and the pairwise distances
between the images $\kappa(\bar x_1)$ (or $\kappa(\bar x_2)$),
$\kappa(\bar x_3)$ and $\kappa(\bar x_4)$) are greater
than  the caliber $\varepsilon_2$ of the approximation.

Type 2. The two pairs $(\kappa(\bar x_1), \kappa(\bar x_2))$ and
$(\kappa(\bar x_3), \kappa(\bar x_4))$ of the images are
$\varepsilon_2$-close in $\RP^s$ and the distance between the
images $\kappa(\bar x_1)$ (or $\kappa(\bar x_2)$) and $\kappa(\bar
x_3)$ (or $\kappa(\bar x_4)$) are greater than
the calibre $\varepsilon_2$ of the approximation. (The
described component is the complement of the regular
$\varepsilon_2$ neighborhood of the triple points manifold of
$d(\RP^s)$.)

Type 3. Images $\kappa(\bar x_1), \kappa(\bar x_2)$ and
$\kappa(\bar x_3)$ on $\RP^s$ are pairwise $\varepsilon_2$-close
in $\RP^s$ and the distance between the images $\kappa(\bar x_1)$
(or $\kappa(\bar x_2)$, or $\kappa(\bar x_3)$) and $\kappa(\bar
x_4)$ is greater than  the caliber
$\varepsilon_2$ of the approximation.

Type 4. All the images $\kappa(\bar x_1), \kappa(\bar x_2)$,
$\kappa(\bar x_3)$ and $\kappa(\bar x_4)$ are pairwise
$\varepsilon_2$-close in $\RP^s$.

Let us prove that there exists a generic $f$ such that the
components of the type 1 and the type 3 are empty. For the
component of the type 3 the proof is analogous to  the proof for
the component of the type 1.

Let us prove that there exists a generic deformation $g_1 \to g_2$
with the caliber $3\varepsilon_3$ such that after this deformation
in the neighborhood $U_{\Delta}^{reg}$ there are no
self-intersection points of $g_2$ obtained by a generic resolution
of triple points of $f$  of the types 1 and 2. Let us start with
the proof for triple points of the type 1.

For a generic small alteration of the immersion $g_2$ inside
$U^{reg}_{\Delta}$ the points of the type 1 of  the triple points
manifold $\Delta_3(f)$ are perturbed into a component of the
self-intersection points on $L^{n-4k}$. This component is
classified by  the following two subtypes:

-- Subtype {\bf a}. Preimages of a point are $(\bar x_2, \bar x_1),
(\bar x_2, \bar x'_1)$.

---Subtype {\bf b}.  Preimages of a point are $(\bar x_1, \bar x'_1),(\bar x_1, \bar x_2)$.

In the formula above the points with the common indeces  have
$\varepsilon_3$-close projections on the corresponding sheet of
$d(\RP^s)$. The two points in a pair form a point on $N^{n-2k}$
and a couple of pairs forms a point on the component of
$L^{n-4k}$.

Let us prove that there exists a $2\varepsilon_3$-small regular
deformation $g_1 \to g_2$, such that the component of $h(L^{n-4k})
\cap U^{reg}_{\Delta}$ of  the subtype {\bf a} is empty. Let
$K^{s-k}$ be the intersection manifold of $f(M^{n-k})$ with
$d(\RP^s)$  (this manifold is immersed into the regular part in $\RP^s$). By a general position
argument, because $2s < n - 2k$, a generic perturbation $r \to r'$ of the immersion
$r: K^{s-k} \looparrowright \RP^s \to \R^n$ is an embedding.
Therefore there exists a $2\varepsilon_2$-small
deformation  of immersed manifold $r(K^{s-k}) \to r'(K^{s-k})$ in $\R^n$, such that
 the regular $\varepsilon_2$-neighborhood of the submanifold $r'(K^{s-k})$ has no self-intersection.
The deformation of the immersed manifolds $r(K^{s-k}) \to
r'(K^{s-k})$ is extended to the deformation of $g_1(N^{n-2k})$ in
the  regular neighborhoods of the constructed one-parameter family
of immersed manifolds. After the described regular deformation the
immersed manifold $g_2(N^{n-2k})$ has no self-intersection
components of the subtype {\bf a}. The case of the
self-intersection of the subtype {\bf b} is analogous.

Let us describe a generic deformation $g_1 \to g_2$ with the
support in $U^{reg}_{\Delta}$ that resolves self-intersection 
corresponding to quadruple points of $f$ of the type 2. This
deformation could be arbitrarily small. After this deformation the
component $\Delta_4(f)$ of the type 2 is resolved into two
components of $L^{n-4k}$ of different subtypes. These two
components will be denoted by $L^{n-4k}_x$, $L^{n-4k}_y$.

The immersed submanifold $g_2(N^{n-2k}) \cap U^{reg}_{\Delta}$ is
divided into two components. The first component is formed by
pairs of points $(\bar x, \bar x')$ with the
$3\varepsilon_3$-close images $(\kappa(\bar x), \kappa(\bar x')$
on $\RP^s$. This component is denoted by $g_2(N^{n-2k}_x)$. The
last component of $g_2(N^{n-2k}) \cap U^{reg}_{\Delta}$ is denoted
by $g_2(N^{n-2k}_y)$. This component is formed by pairs of points
$(\bar x, \bar x')$ with the projections $(\kappa(\bar x),
\kappa(\bar x'))$ on different sheets of $\RP^s$.

The component  $L^{n-4k}_{x\downarrow}$ is defined by pairs $(\bar
x_1, \bar x'_1),(\bar x_2, \bar x'_2)$. The component
$L^{n-4k}_{y}$ is defined by pairs $(\bar x_1, \bar x_2),(\bar
x'_1, \bar x'_2)$. A common index  of points in the pair means
that the images of the points are $\varepsilon_3$-close on
$\RP^s$.  Each  pair determines a point on $N^{n-2k}$ with the
same image of $g_2$. It is easy to see that the component
$L^{n-4k}_{x\downarrow}$ is the self-intersection of
$g_2(N^{n-2k}_x)$ and the component $L^{n-4k}_{y}$ is the
self-intersection of $g_2(N^{n-2k}_y)$.

It is easy to see that the structure groups of the components 
agree with the corresponding  subgroup described in the lemma. The
component $L^{n-4k}_{x \downarrow}$  admits a reduction of the
structure group to the subgroup $\I_{2,x \downarrow} \subset \Z/2
\int \D_4$. The component $L^{n-4k}_{y}$ admits a reduction of the
structure group to the subgroup $\I_{2,y}$. Moreover, it is easy
to see that the covering $\tilde L^{n-4k}_{x \downarrow}$ over
$L^{n-4k}_x$ induced by the epimorphism $\omega: \Z/2 \int \D_4
\to \Z/2$ with the kernel $\I_{2,x} \subset \Z/2 \int \D_4$ is
naturally diffeomorphic to $L^{n-4k}_y$. Also it is easy to see that
this  diffeomorphism  agrees with the transformation $OP$ of the
structure groups of the framing over the components.

The last component of $L^{n-4k}$ is immersed in the
$\varepsilon_2$-neighborhood of $d(\RP^s)$ outside of
$U^{reg}_{\Delta}$ and will be denoted by $L^{n-4k}_z$. The
structure group of the framing of this component is $\I_{2,z}$.
Lemma 1 is proved.

\subsubsection*{The last part of the proof of the Theorem 1}

Let us construct a pair of polyhedra $(P',Q') \subset \R^n$,
$dim(P') =2s-n=n-2k-q-2$, $dim(Q') = dim(P')-1$. Obviously,
$dim(P')< 2k -1$. Take a generic mapping $d': \RP^s \to \R^n$. Let
us consider the submanifold with boundary $(\Delta'^{reg},
\partial \Delta'^{reg}) \subset \R^n$ (see the denotation in 
Lemma 1). Let $\eta_{\Delta'^{reg}}: (\Delta'^{reg},\partial
\Delta'^{reg}) \to (K(\D_4,1),K(\I_b,1))$ be the classifying
mapping for the double point self-intersection manifold of $d$.

 By a standard argument we may modify the mapping $d$ into $d'$ such that the
mapping $\eta_{\Delta^{reg}}$ is a homotopy equivalence of 
 pairs up to the dimension $q+1$. After this
modification $d' \to d$ we define
 $(P,Q)=(\Delta^{reg},\partial \Delta^{reg}) \subset \R^n$ and the
 mapping $\eta_{\Delta^{reg}}$ is a $(q+1)$-homotopy equivalence.

The subpolyhedron $Q$  is equipped with two cohomology classes
$\kappa_{Q,1}, \kappa_{Q,2} \in H^1(Q;\Z/2)$.  Because $\Sigma$ is
a submanifold in $\RP^s$, the restriction of the characteristic
class $\kappa \in H^1(\RP^s;\Z/2)$ to $H^1(\Sigma;\Z/2)$ is
well-defined. The inclusion $i_Q: Q \subset U_{\Sigma}$ determines
the cohomology class $(i_Q)^{\ast}(\kappa) \in H^1(Q;\Z/2)$.
The cohomology class $\kappa_{Q,1}$ is defined as the
characteristic class of the canonical double points covering over
$\Sigma$. The class $\kappa_{Q,2}$ is defined by the formula
$\kappa_{Q,2} = (i_Q)^{\ast}(\kappa) + \kappa_{Q,1}$.

The immersed manifold (with boundary) $(N^{n-2k} \cap U_{\Sigma})
\looparrowright U_{\Sigma}$ is equipped with an $\I_b$-framing.
Obviously the classes $\kappa_{Q,1}, \kappa_{Q,2} \in
H^1(U_{\Sigma};\Z/2)=H^1(Q;\Z/2)$ restricted to
$H^1(g_2(N^{n-2k}_{ext});\Z/2)$ ( recall that
$g_2(N^{n-2k}_{ext})=g_2(N^{n-2k}) \cap (\R^n \setminus
U_{\Delta})$)   agree with the two generated cohomology
classes  $\rho_1, \rho_2$ of the $\I_b$-framing correspondingly.

 Let us define the immersion $g: N^{n-2k} \looparrowright \R^n$ with 
 $\I_b$-control over $(P,Q)$. Let us start with the immersion $g_2: N^{n-2k} \looparrowright \R^n$
 constructed in the lemma. By a $2\varepsilon_2$--small generic
 regular deformation we may deform the immersion $g_2$ into $g_3$,
 such that this deformation pushes the component $g_2(N^{n-2k}_x)$ out of $U^{reg}_{\Delta}$.
 Therefore the component $L^{n-4k}_{x \downarrow} \subset L^{n-4k}$
 of the self-intersection of $g_2$ is also deformed out of
 $U^{reg}_{\Delta}$.

The immersed manifold (with boundary) $g_3(N^{n-2k}) \cap (\R^n
\setminus U^{reg}_{\Delta})$ is equipped with an $\I_b$-framing
of the normal bundle.
 Obviously,  the classes
$\kappa_{Q,1}, \kappa_{Q,2} \in H^1(U_{\Sigma};\Z/2)=H^1(Q;\Z/2)$,
restricted to $H^1(g_2(N^{n-2k}) \cap U_{\Delta};\Z/2)$,  agree
with the two generated cohomological classes of the
$\I_b$-framing. The immersed manifold $g_3(N^{n-2k}) \cap
U^{reg}_{\Delta}$ coincides with $g_2(N^{n-2k}_{y})$ and has the
general structure group of the framing. This immersed manifold has
the self-intersection manifold (with boundary) $h(L^{n-4k}) \cap
U^{reg}_{\Delta}$ with the reduction of the structure group to the
pair of the subgroups $(\I_{2,y},\I_3)$.

Let us prove that the immersed manifold (with boundary)
$h(L^{n-4k}) \cap U^{reg}_{\Delta}$ is $\Z/2 \int \D_4$-framed
cobordant (relative to the boundary) to a $\Z/2 \int \D_4$-framed
manifold decomposed into the disjoint union of a closed $\Z/2 \int
\D_4$-framed manifold that is the image of the transfer
homomorphism $\omega^!$ and a relative $\I_3$-framed manifold.

Take a $\Z/2 \int \D_4$-framed manifold $(\tilde L^{n-4k}, \tilde
\Psi, \tilde \zeta)$ that is defined as the image of $\Z/2 \int
\D_4$-framed manifold $(L^{n-4k}, \Psi, \zeta)$ by the transfer
homomorphism (a double covering) with respect to the cohomology
class $\omega \in H^1(\Z/2 \int \D_4;\Z/2)$. Recall that the
manifold $\tilde L^{n-4k}$ is obtained by  gluing the manifold
$\tilde L^{n-4k}_{x} \cup \tilde L^{n-4k}_y$ with the manifold
$\tilde L^{n-4k}_z$ along the common boundary $\tilde
\Lambda^{n-4k-1}$. Note that the group of the framing of the last
manifold $\tilde \Lambda^{n-4k-1}_z$ is the subgroup $\I_3 \subset
\Z/2 \int \D_4$.

Let $OP\alpha$ be the $\Z/2 \int \D_4$--framed immersion  obtained
from an arbitrary $\Z/2 \int \D_4$-framed immersion $\alpha$ by
changing  the structure group of the framing by the
transformation $OP$. The $\Z/2 \int \D_4$-framed manifold (with
boundary) $(\tilde L^{n-4k}_y, \tilde \Psi_y, \tilde \zeta_y)$
coincides with  the two disjoint copies of $\Z/2 \int \D_4$-framed
manifold (with boundary) $OP(\tilde L^{n-4k}_y, \tilde \Psi_y,
\tilde \zeta_y)$.

Let us put $\alpha_1=-OP(\tilde L^{n-4k}, \tilde \Psi, \tilde
\zeta)$. Let us define the sequence of  $\Z/2 \int \D_4$-framed
immersions $\alpha_2 = -2 OP\alpha_1$, $\alpha_3 = -2 OP\alpha_2$,
$\dots$, $\alpha_j = -2 OP\alpha_{j-1}$.

Obviously,  the $\D/4 \int \Z/2$-framed immersion $\alpha_1 +
\alpha_2 = \alpha_1 + 2OP\alpha_1^{-1}$ is represented by 3 copies
of the manifold $\tilde L^{n-4k}$. The second and the third copies
are obtained from the first copy by the mirror image  and the
changing of structure group of the framing. The manifold
$-OP[\tilde L^{n-4k}] \cup 2[\tilde L^{n-4k}]$ contains, in
particular, a copy of $-OP[\tilde L^{n-4k}_x]$ inside the first
component and the union $[\tilde L^{n-4k}_y \cup L^{n-4k}_y]$ of
the mirror two copies of $-OP[\tilde L^{n-4k}_x]$  in the second
and the third component. Therefore the manifold $-OP[\tilde
L^{n-4k}] \cup 2[ \tilde L^{n-4k}]$ is $\Z/2 \int \D_4$-framed
cobordant to a $\Z/2 \int \D_4$-framed manifold, obtained by 
gluing the union of a copy of $-OP[\tilde L^{n-4k}_x]$ and 4
copies of $\tilde L^{n-4k}_y$ by a $\I_3$-framing manifold along
the boundary. This cobordism is relative with respect to the
submanifold $-OP[\tilde L^{n-4k}_z] \cup 2[L^{n-4k}_z] \subset
-OP[L^{n-4k}] \cup 2L^{n-4k}$.

By an analogous argument it is easy to prove that  the element
$\aleph= \sum_{j=1}^{j_0} \alpha_j$ is $\Z/2 \int \D_4$-framed
cobordant to the manifold obtained by gluing the union
$-OP[\tilde L^{n-4k}_x] \cup 2^j(-OP)^{j-1}[\tilde L^{n-4k}_y]$ by
an $\I_3$-manifold along the boundary. Moreover, this cobordism is
relative with respect to all copies of $\tilde L^{n-4k}_z$
(with various orientations). If $j_0$ is great enough, the
manifold (with $\I_3$-framed boundary) $2^j(-OP)^{j_0-1}[\tilde
L^{n-4k}_y]$ is cobordant relative to the boundary to an
$\I_3$-framed manifold.

Therefore the manifold $L^{n-4k}_y$ is $\Z/2 \int \D_4$-framed
cobordant relative to the boundary to the union of an $\I_3$-framed
manifold with the same boundary and a closed manifold that is the
double cover with respect to $\omega$ over a $\Z/2 \int
\D_4$-framed manifold. This cobordism is realized as a cobordism
of the self-intersection  of  a $\D_4$-framed immersion with 
support inside $U^{reg}_{\Delta}$. This cobordism joins the
immersion $g_3$ with a  $\D_4$--framed immersion $g_4$. After an
additional deformation of $g_4$ inside a larger neighborhood of
$\Delta^{reg}$ the relative $\I_b$-submanifold of the
self-intersection manifold of $g_4$ is deformed outside of
$U^{reg}_{\Delta}$. The $\D_4$-framed immersion obtained as the
result of this cobordism admits an $\I_b$-control. The Theorem 1
is proved.

\section{ An $\I_4$-structure
(a cyclic structure) of a $\D_4$-framed immersion}

Let us describe the subgroup $\I_4 \subset \Z/2 \int \D_4$. This
subgroup is isomorphic to the group $\Z/2 \oplus \Z/4$. Let us
recall that the group $ \Z/2 \int \D_4$ is the transformation
group of $\R^4$ that permutes the $4$-tuple of the coordinate lines
and two planes $(f_1, f_2)$, $(f_3, f_4)$ spanned by the vectors
of the standard base $(f_1, f_2, f_3, f_4)$ (the planes can 
remin fixed or be permuted by a transformation).

Let us denote the generators of $\Z/2 \oplus \Z/4$ by $l$, $r$
correspondingly. Let us describe the transformations of $\R^4$
given by each generator. Consider a new base $(e_1, e_2,
e_3, e_4)$, given by $e_1=f_1+f_2$, $e_2=f_1-f_2$, $e_3=f_3+f_4$,
$e_4=f_3-f_4$. The generator $r$ of order 4 is represented by
the rotation in the plane $(e_2,e_4)$ through the angle $\frac{\pi}{2}$
and the reflection in the plane $(e_1,e_3)$ with respect to the
line $e_1+e_3$. The generator $l$ of order 2 is represented by
the central symmetry  in the plane $(e_1, e_3)$.

Obviously, the described representation of $\I_4$ admits 
invariant (1,1,2)-dimensional subspaces. We will denote 
subspaces by  $\lambda_1, \lambda_2, \tau$.

The lines  $\lambda_1, \lambda_2$ are generated by the
vectors $e_1+e_3$, $e_1-e_3$ correspondingly. The subspace $\tau$
is generated by the vectors $e_2, e_4$. The generator $r$ acts by
the reflection in $\lambda_2$ and by the rotation in $\tau$ throught the
angle $\frac{\pi}{2}$. The generator $l$ acts by reflections in
the subspaces  $\lambda_1$, $\lambda_2$.

In particular, if the structure group $\Z/2 \int \D_4$ of a
4-dimensional bundle $\zeta: E(\zeta) \to L$ admits a reduction
to the subgroup $\I_4$, then the bundle is decomposed into the
direct sum $\zeta = \lambda_1 \oplus \lambda_2 \oplus \tau$ of
$1,1,2$--dimensional subbundles.

\subsubsection*{Definition 6}

Let $(g: N^{n-2k} \looparrowright \R^n, \Xi_N, \eta)$ be an
arbitrary $\D_4$-framed immersion. We shall say that this
immersion is an $\I_b$--immersion (or a cyclic immersion), if the
structure group $\Z/2 \int \D_4$ of the normal bundle over the
double points manifold $L^{n-4k}$ of this immersion admits a
reduction to the subgroup $\I_4 \subset \Z/2 \int \D_4$. In this
definition we assume that the pairs $(f_1, f_2)$, $(f_3,f_4)$ are
the vectors of the framing for the two sheets of the
self-intersection manifold at a point in the double point manifold
$L^{n-4k}$.
\[  \]

In particular, for a cyclic $\Z/2 \int \D_4$-framed immersion
there exists the mappings $\kappa_a: L^{n-4k} \to K(\Z/2,1)$,
$\mu_a: L^{n-4k} \to K(\Z/4,1)$ such that the characteristic
mapping $\zeta: L^{n-4k} \to K(\Z/2 \int \D/4,1)$ of the $\Z/2
\int \D_4$-framing of the normal bundle over $L^{n-4k}$ is reduced
to a mapping with the target $K(\I_b,1)$ such that the following
equation holds: 
$$\zeta = i(\kappa_a \oplus \mu_a),$$ 
where  $i: \Z/2
\oplus \Z/4 \to \I_4$ is the prescribed isomorphism.

The following Proposition is proved by a straightforward
calculation.

\subsubsection*{Proposition 2}
Let $(g,\Psi_N, \eta)$ be a $\D_4$--framed immersion, that is a
cyclic immersion. Then the Kervaire invariant, appearing as the top line of the
diagram (7), can be calculated by following formula:
$$\Theta_a=<\kappa_a^{\frac{n-4k}{2}}\mu_a^{\ast}(\tau)^{\frac{n-4k-2}{4}}
\mu_a^{\ast}(\rho);[L]>, \eqno(8)$$ where $\tau \in
H^2(\Z/4;\Z/2)$, $\rho \in H^1(\Z/4;\Z/2)$  are the generators.
\[  \]

\subsubsection*{Proof of Proposition 2}
 Let us consider the subgroup of  index 2,
 $\I_b \subset \I_4$.  This subgroup is the kernel of the epimorphism $\chi':
\I_4 \to \Z/2$, that is the restriction of the characteristic
class $\chi: \Z/2 \int \D_4 \to \Z/2$ of the canonical double
cover $\bar L \to L$ to the subgroup $\I_4 \subset \Z/2 \int
\D_4$. Obviously, the characteristic number (8) is calculated by
the formula
$$ \Theta_a=<\hat \kappa_a^{\frac{n-4k}{2}} \hat
\rho_a^{\frac{n-4k}{2}};\bar L>, \eqno(9)$$ where the
characteristic class $\hat \kappa_a \in H^1(\bar L;\Z/2)$ is
induced from the class $\kappa_a \in H^1(L;\Z/2)$ by the canonical
cover $\bar L \to L$, and the class $\hat \rho_a \in H^1(\bar
L;\Z/2)$ is obtained by the transfer of the class $\rho \in
H^1(L;\Z/4)$.

Note that $\hat \kappa_a = \tau_1$, $\hat \rho_a = \tau_2$, where
$\tau_1$, $\tau_2$ are the two generating $\I_b$--characteristic
classes. Therefore $\hat \kappa_a \hat \rho_a = \tau_1 \tau_2 =
w_2(\eta)$, where $\eta$ is the two-dimensional bundle that
determines the $\D_4$--framing (over the submanifold $\bar
L^{n-4k} \subset N^{n-2k}$ this framing admits a reduction to an
$\I_b$-framing) of the normal bundle for the immersion $g$ of 
$N^{n-2k}$ into $\R^n$.

Therefore the characteristic number, given by the formula (8) in
the case when the $\Z/2 \int \D_4$ framing over $L^{n-4k}$ is
reduced to an $\I_4$-framing, coincides with the characteristic
number, given by the formula (9). Proposition 2 is proved.

\subsubsection*{Definition 7}
We shall say that a $\D_4$-framed immersion $(g,\Xi_N,\eta)$
admits a $\I_4$--structure (a cyclic structure), if for the double points
manifold $L^{n-4k}$ of $g$ there exist mappings $\kappa_a:
L^{n-4k} \to K(\Z/2,1)$, $\mu_a: L^{n-4k} \to K(\Z/4,1)$ such that
the characteristic number (8) coincides with  Kervaire
invariant, see Definition 2.
\[  \]

\subsubsection*{Theorem 2}

Let  $(g, \Psi, \eta)$ be a $\D_4$-framed immersion, $g: N^{n-2k}
\looparrowright \R^n$, that represents a regular cobordism class
in the image of the homomorphism $\delta: Imm^{sf}(n-k,k) \to
Imm^{\D_4}(n-2k,2k)$, $n-4k=62$, $n=2^l-2$, $l \ge 13$, and assume the
conditions of the Theorem 1 hold, i.e. the residue class
$\delta^{-1}(Imm^{sf}(n-k,k)$ (this class is defined modulo odd
torsion) contains a skew-framed immersion that admits a
retraction of  order  $62$.

Then in the $\D_4$-framed cobordism class $[(g, \Psi, \eta)] =
\delta[(f, \Xi, \kappa)] \in Imm^{\D_4}(n-2k,2k)$ there exists a
$\D_4$-framed immersion that admits an $\I_4$--structure (a cyclic structure).
\[ \]

\section{Proof of  Theorem 2}
 Let us formulate the Geometrical Control Principle for $\I_b$--controlled immersions.

Let us take an $\I_b$--controlled immersion (see Definition 4)
$(g,\Xi_N,\eta;(P,Q),\kappa_{Q,1}, \kappa_{Q,2})$, where $g:N
\looparrowright \R^n$ is a $\D_4$-framed immersion, equipped with a control mapping
over a polyhedron  $i_P: P \subset \R^n$, $dim(P)=2k-1$; $Q
\subset P$  $dim(Q)=dim(P)-1$. The characteristic classes $\kappa_{Q,i} \in H^1(Q;\Z/2)$,
$i=1,2$  coincide with characteristic classes
$\kappa_{i,N_Q} \in N_Q^{n-2k-1}$ by means of the mapping $\partial N^{n-2k}_{int} = N^{n-2k}_Q \to Q$,
where $N^{n-2k}_{int} \subset N^{n-2k}$,
$N^{n-2k}_{int}=g^{-1}(U_P)$, $U_P \subset \R^n$.

\subsubsection*{Proposition 3. Geometrical Control Principle for $\I_b$--controlled immersions}

Let $j_P : P \subset \R^n$ be an arbitrary embedding; such an
embedding is unique up to isotopy by a dimensional reason,
because $2dim(P)+1=4k-1<n$. Let $g_1: N^{n-2k} \to \R^n$ be an
arbitrary mapping, such that the restriction $g_1 \vert_{N_{int}}:
(N^{n-2k}_{int}, N^{n-2k-1}_Q) \looparrowright (U_P, \partial
U_P)$ is an immersion (the restriction $g \vert_{N^{n-2k-1}_Q}$ is
an embedding) that corresponds to the immersion $g
\vert_{N^{n-2k}_{int}}: (N^{n-2k}_{int},N^{n-2k-1}_Q)
\looparrowright (U_P, \partial U_P)$ by means of the standard
diffeomorphism of the regular neighborhoods $U_{i_P}=U_{j_P}$ of
subpolyhedra $i(P)$ and $j(P)$. (For a dimension reason there is
a standard diffeomorphism of $U_{i_P}$ and $U_{j_P}$ up to an
isotopy.)

Then for an arbitrary  $\varepsilon
>0$ there exists an immersion $g_{\varepsilon}: N^{n-2k} \looparrowright
\R^n$ such that $dist_{C^0}(g_1,g_{\varepsilon})<\varepsilon$ and such that
 $g_{\varepsilon}$ is regular homotopy to an immersion  $g$ and the restrictions
 $g_{\varepsilon} \vert _{N^{n-2k}_{int}}$ and $g_1 \vert_{N^{n-2k}_{int}}$  coincide.
\[  \]

We start the proof of Theorem 2 with the following construction. Let us consider the
manifold $Z=S^{\frac{n}{2}+64}/i \times \RP^{\frac{n}{2}+64}$.
This manifold is the direct product of the standard lens
space $(mod 4)$ and the projective space. The cover $p_Z: \hat Z \to Z$ over this
manifold with the covering space $\hat Z = \RP^{\frac{n}{2}+64}
\times \RP^{\frac{n}{2}+64}$ is well-defined.

Let us consider in the manifold $Z$ a family of submanifolds
$X_i$, $i=0, \dots, \frac{n+2}{64}$ of the codimension
$\frac{n+2}{2}$, defined by the formulas $X_0 =
S^{\frac{n}{2}+64}/i \times \RP^{63}$, $X_1= S^{\frac{n}{2}+32}/i
\times \RP^{95}, \dots$, $X_j = S^{\frac{n}{2} - 32(j-2)-1}/i
\times \RP^{32(j+2)-1}, \dots$, $X_{\frac{n+2}{64}} = S^{63}/i
\times \RP^{\frac{n}{2}+64}$. The embedding of the corresponding
manifold in $Z$ is defined by the Cartesian product of the two
standard embeddings.

The union of the submanifolds
 $\{X_i\}$ is a stratified submanifold (with singularities)
 $X \subset Z$ of the dimension
$\frac{n}{2}+127$, the codimension of maximal singular strata in
$X$ is equal to $64$. The covering  $p_X: \hat X \to X$, induced
from the covering $p_Z: \hat Z \to Z$ by the inclusion $X \subset
Z$, is well-defined. The covering space $\hat X$ is a stratified
manifold (with singularities) and decomposes into the union of the
submanifolds $\hat X_0 = \RP^{\frac{n}{2}+64} \times \RP^{63},
\dots, \hat X_j = \RP^{\frac{n}{2} - 32(j-2)} \times \RP^{32(j+2)
-1}, \dots, \hat X_{\frac{n+2}{64}}=\RP^{63} \times
\RP^{\frac{n}{2}+64}$.
 Each
manifold $\hat X_i$ of the family is the $2$-sheeted covering
space over the manifold $X_i$ over the first coordinate. Let us
define $d_1(j)= \frac{n}{2} - 32(j-2)$, $d_2(j)=32(j+2)-1$. Then
the formula for $X_i$ is the following:  $X_j=\RP^{d_1(j)} \times
\RP^{d_2(j)}$.

The cohomology classes  $\rho_{X,1} \in H^1(X;\Z/4)$,
$\kappa_{X,2}\in H^1(X;\Z/2)$ are well-defined. These classes are
induced from the generators  of the groups $H^1(Z;\Z/4)$,
$H^1(Z;\Z/2)$. Analogously, the cohomology classes
 $\kappa_{\hat
X,i}\in H^1(\hat X;\Z/4)$, $i=1,2$ are well-defined. The cohomology class $\kappa_{\hat X,1}$
is induced from the class $\rho_{X,1} \in H^1(X;\Z/4)$ my means of the transfer homomorphim,
and $\kappa_{\hat X,2} = (p_X)^{\ast}(\kappa_{X,2})$.

Let us define for an arbitrary $j= 0, \dots, (\frac{n+2}{64})$ the
space $J_j$ and the mapping $\varphi_j: X_j \to J_j$. We denote by $Y_1(k)$ the space $S^{31}/i \ast \dots
\ast S^{31}/i$ of the join of $k$ copies,  $k=1, \dots
,(\frac{n+2}{64}+1)$, of the standard lens space $S^{31}/i$. Let
us denote by $Y_2(k)$, $k=2, \dots, (\frac{n+2}{64}+2)$, $Y_2(k) =
\RP^{31} \ast \dots \ast \RP^{31}$ the joins of the $k$ copies of
the standard projective space $\RP^{31}$. Let us define $J_j
=Y_1(\frac{n+2}{64}-j+2)) \times Y_2(j+2)$
 $Q =
Y_1(\frac{n+2}{64}+2) \times Y_2(\frac{n+2}{64}+2)$. For a given
$j$ the natural inclusions $J_j \subset Q$ are well-defined. Let
us denote the union of the considered inclusions by $J$.

The mapping
$\varphi_j: X_j \to J_j$ is well-defined as the Cartesian product
 of the two following mappings.
 On the first coordinate the mapping is defined as the composition of the standard 2-sheeted covering
$\RP^{d_1(j)} \to S^{\frac{n}{2}-64(j-1)}/i$ and the natural
projection $S^{d_1(j)}/i \to Y_1(d_1(j))$. On the second
coordinate the mapping is defined by the natural projection
$\RP^{d_2(j)} \to Y_2(j+1)$.

The family of mappings
 $\varphi_j$ determines the mapping
$\varphi: \hat X \to J$, because the restrictions of any two
mappings to the common subspace in the origin coincide.

For $n+2 \ge 2^{13}$ the space  $J$ embeddable into the
Euclidean $n$-space by an embedding $i_J: J \subset \R^n$. Each
space $Y_1(k)$, $Y_2(k)$ in the family is embeddable into the
Euclidean  $(2^6k -1-k)$--space. Therefore for an arbitrary $j$
the space $J_j$ is embaddable into the Euclidean space of 
dimension $n+126-\frac{n+2}{64}$. In particular, if $n+2 \ge
2^{13}$ the space $J_j$ is embeddable into $\R^n$. The image of an
arbitrary intersection of the two embeddings in the family belongs
to the standard coordinate subspace. Therefore the required
embedding $i_J$ is defined by the gluing of embeddings in the
family.

Let us describe the mapping $\hat h: \hat X \to \R^n$. By
$\varepsilon$ we denote the radius of a (stratified) regular
neighborhood of the subpolyhedron  $i_J(J) \subset \R^n$. Let us
consider a small positive $\varepsilon_1$, $\varepsilon_1 <<
\varepsilon$, (this constant will be defined below in the proof of
Lemma 4) and let us consider a generic $PL$
$\varepsilon_1$--deformation of the mapping $i_J \circ \varphi:
\hat X \to J \subset \R^n$. The result of the deformation is
denoted by $\hat h: \hat X \to \R^n$.

Let us define the positive integer  $k$ from the equation
$n-4k=62$. In the prescribed regular homotopy class of an
$\I_b$-controlled immersion $f: N^{n-2k} \looparrowright \R^n$ we
will construct another $\I_b$--controlled immersion $g: N^{n-2k}
\looparrowright \R^n$  that admits a $\I_b$--structure.

Let the immersion $f$ be controlled over the embedded
subpolyhedron $\psi_P: P \subset \R^n$. Let $\psi_Q: Q \to \hat X$
be a generic mapping such that $\kappa_{Q,i} = \psi_Q \circ
\kappa_{\hat X,i}$, $i=1,2$. By the previous definition the
manifolds $N^{n-2k}_{int}$, $N^{n-2k}_{ext}$ with the common
boundary $N^{n-2k-1}_Q$, $N^{n-2k} = N^{n-2k}_{int}
\cup_{N^{n-2k-1}_Q} N^{n-2k}_{ext}$ are well-defined.

Let $\eta: N^{n-2k}_{ext} \to K(\I_b,1) \subset K(\D_4,1)$ be the
characteristic mapping of the framing $\Xi_N$, restricted to
$N^{n-2k}_{ext} \subset N^{n-2k}$. The restriction of this mapping
to the boundary $\partial N^{n-2k}_{ext}= N^{n-2k-1}_Q$ is given
by the composition $\partial N^{n-2k-1}_Q \to Q \to K(\I_b,1)
\subset K(\D_4,1)$. The target space for the mapping $\eta$ is the
subspace $K(\I_b,1) \subset K(\D_4,1)$. This mapping is determined
by the cohomology classes $\kappa_{N^{n-2k}_{ext},s} \in
H^1(N^{n-2k}_{ext},Q;\Z/2)$, $s=1,2$.

Let us define the mapping $\lambda: N^{n-2k}_{ext} \to \hat X$ by
the following conditions. This mapping transforms the
cohomology classes  $\kappa_{\hat X,i}$ into the classes
$\kappa_i \in H^1(N^{n-2k}_{ext};\Z/2)$ and also the restriction
$\lambda \vert_{N^{n-2k-1}_Q}$ coincides with the composition of
the projection $N^{n-2k-1}_{Q} \to Q$ and the mapping $\psi_Q: Q
\to \hat X$. The boundary conditions for the mapping $\psi_Q$ are
$\kappa_{Q,i} = \psi_Q \circ \kappa_{\hat X,i}$, $i=1,2$. The
submanifold with singularities $\hat X \subset \hat Z$ contains
the skeleton of the space $\hat Z$ of the dimension
$\frac{n}{2}+62$. Because $n-2k=\frac{n}{2}+31$, the mapping
$\lambda$ is well-defined.

Let us denote the composition
 $\hat h \circ \lambda: N^{n-2k}_{ext} \to \hat X \to
\R^n$ by $g_1$. Let us denote the mapping
$ \hat h \circ \psi_Q: Q \to \hat X \to \R^n$ by $
\varphi_Q$. One can assume that the mapping  $\varphi_Q$ is an embedding.
Moreover, without loss of  generality one may assume that this embedding
is extended to a generic embedding  $\varphi_P: P \subset
\R^n$ such that the embedded polyhedron $\varphi_P: P \subset
\R^n$ does not intersect  $g_1(N^{n-2k}_{ext})$.

Let us denote by $U_{\varphi}(P)$ a regular neighborhood of the
subpolyhedron  $\varphi_P(P) \subset \R^n$ (we may assume that the
radius of this neighborhood is equal to $\varepsilon$). Up to an
isotopy a regular neighborhood $U_{\varphi}(P)$ is well-defined,
in particular, this neighborhood does not depend on the choice of
a regular embedding of $P$, moreover $U_{\varphi}(P)$ and $U(P)$ are
diffeomorphic.

Without loss of  generality after an additional small
deformation we may assume that the restriction $g_1
\vert_{N^{n-2k}_{int}}$ is a regular immersion $g_1:
N^{n-2k}_{int} \subset \R^n$ with the image inside
$U_{\varphi}(P)$. In particular, the restriction of $g_1$ to the
boundary $N^{n-2k-1}_Q = \partial(N^{n-2k}_{int})$ is a regular
embedding $N^{n-2k-1}_Q \subset \partial U(P)$. The immersion $g_1
\vert_{N_{int}}$ is conjugated to the immersion  $f
\vert_{N_{int}}$ by means of a diffeomorphism  of $U_{\varphi}(P)$
with $U(P)$.

By Proposition 3, for an arbitrary $\varepsilon_2 >0$,
$\varepsilon_2 << \varepsilon_1 << \varepsilon$, there exists an
immersion $g : N^{n-2k} \looparrowright \R^n$ in the regular
homotopy class of  $f$, such that $g$ coincides with $g'$ (and
with  $g_1$) on $N^{n-2k}_{int}$ and, moreover, $dist(g,g_1)<
\varepsilon_2$.

Let us consider the self-intersection manifold $L^{n-4k}$ of the
immersion $g$. This manifold is a submanifold in $\R^n$. Let us
construct  the mappings $\kappa_a: L^{n-4k}
\to K(\Z/2,1)$, $\mu_a: L^{n-4k} \to K(\Z/4,1)$. Then we check the
conditions (8) and (9).

The manifold $L^{n-4k}$ is naturally divided into two components.
The first component
  $L^{n-4k}_{int}$ is inside $U_{\varphi_P}(P)$. The last component
  (we will denote this component again by $L^{n-4k}$) consists of the last self-intersection points.
  This component is outside the $\varepsilon$--neighborhood
  of the submanifold with singularities $h(X)$.
The mappings $\kappa_a$, $\mu_a$ over $L^{n-4k}_{int}$ are defined as the trivial mappings.
Let us define the mappings $\kappa_a$, $\mu_a$ on $L^{n-4k}$.

Let us consider the mapping
 $\varphi: \hat X \to J$ and the singular set (polyhedron) $\Sigma$ of this mapping.
This is the subpolyhedron
 $ \Sigma
\subset \{ \hat X^{(2)}=  \hat X \times \hat X \setminus
\Delta_{\hat X} / T'\}$, where $T': \hat X^({2}) \to \hat
X^{(2)}$-- is the involution of coordinates in the delated product
$\hat X^{(2)}$ of the space $\hat X$. The subpolyhedron (it is
convenient to view this polyhedron as a manifold with
singularities) $\Sigma$ is naturally decomposed into the union of
the subpolyhedra $\Sigma(j)$, $j=0, \dots, \frac{n+2}{128}$. The
subpolyhedron $\Sigma(j)$ is the singular set of the mapping
$\varphi(j): \RP^{d_1(j)} \times \RP^{d_2(j)} \to S^{d_1(j)}/i
\times \RP^{d_2(j)} \to J_j$.  This subpolyhedron consists of the
singular points of the mapping $\varphi$ in the inverse image
$(\varphi)^{-1}(J_j)= \RP^{d_1(j)} \times \RP^{d_2(j)}$ of the
subspace $J_j \subset J$.

Let us consider the subspace $\Sigma^{reg} \subset \Sigma$,
consisting of points on strata of length 0 (regular strata) and
of length 1 (singular strata of the codimension 32) after the
regular  $\varepsilon_2$ --neighborhoods ($\varepsilon_2 <<
\varepsilon_1$) of the diagonal $\Delta^{diag}$ and the
antidiagonal $\Delta^{antidiag}$  of $\Sigma^{reg}$ are
cut out.

The manifold with singularities  $\Sigma^{reg}$ admits a natural
compactification (closure) in the neighborhood of
  $\Delta^{diag}$ and $\Delta^{antidiag}$; the result of the compactification will be denoted by $K_{reg}$.

The space $RK$, called the space of  resolution of singularities,
equipped with the natural projection $RK \to K_{reg}$ is defined
by the analogous construction; see the short English translation of [A1], Lemma 7. The
cohomology classes $\rho_{RK,1} \in H^1(RK;\Z/4)$,
$\kappa_{RK,2} \in H^1(RK;\Z/2)$ are well-defined. The
cohomology classes $\kappa_{K_{reg},1} \in H^1(RK;\Z/2)$,
$\kappa_{RK,1} \in H^1(RK;\Z/2)$ are the images of the class
$\kappa_{\Sigma,1} \in H^1(\Sigma;\Z/2)$ with respect to the
inclusion $K_{reg} \subset \Sigma$ and the projection $RK \to
K_{reg}$. The class classifies the transposition of the two
non-ordered preimages of a point in the singular set.

Let us consider the restrictions of the classes $\kappa_{K_{reg},1},
\kappa_{RK,1}, \kappa_{\Sigma,1}$ to  neighborhoods of the
diagonal and the antidiagonal. The natural projection
$\Delta^{diag} \to \hat X$ is well-defined. The restrictions of the
classes $\rho_1$ and $\kappa_2$ to  neighborhoods of the
diagonal coincide with the restrictions of the classes
$\rho_{\hat X,1} \in H^1(\hat X;\Z/4)$, $\kappa_{\hat X,2} \in
H^1(\hat X;\Z/2)$. (These classes $\rho_{\hat X,1}, \kappa_{\hat X,2}$
are extended to 
neighborhoods of the diagonal).

Let us recall that the mapping $\hat h: \hat X \to \R^n$ is
defined as the result of an $\varepsilon_1$--small regular
deformation of the mapping $\hat X \to
 X \stackrel{h}{\longrightarrow}
 \R^n$. The singular set of the mapping $\hat h$ will be denoted
by $\Sigma_{\hat h}$. This is a $128$--dimensional polyhedron, or
a manifold with singularities in the codimensions
$32,64,96,128$. Moreover, the inclusion $\Sigma_{\hat h} \subset
\hat X^{(2)}$ is well-defined. The image of this inclusion is in
the regular $\varepsilon_1$--small neighborhood of the singular
polyhedron $\Sigma \subset X^{(2)}$.

Let us denote by $\Sigma_{\hat h}^{reg}$ the part of the singular
set after cutting out the regular
$\varepsilon_1$--neighborhood of the points in singular strata of
length at least 2 (of the codimension 64) and
self-intersection points of all singular strata (these strata are
also of the codimension 64). The boundary $\partial \Sigma_{\hat
h}$ is a submanifold with singularities in $\hat X$ and therefore. by a general position argument,
we may also assume that the boundary $\partial \Sigma_{\hat
h}^{reg}$ is a regular submanifold with singularities in $\hat X$.

Additionally, by  general position arguments, the intersection of
the image
 $Im(\lambda
(N^{n-2k}_{ext}))$ inside the singular set $\Sigma_{\hat h}$ (this
is a polyhedron of the dimension $62$) on $X$ are outside (with
respect to the caliber $\varepsilon$) of the projection  of the
singular submanifold
 with
singularities (this singular part is of the codimension $64$) in
the complement of the regular submanifold with singularities
$\Sigma_{\hat h}^{reg} \subset \Sigma_{\hat h}$. Therefore the
image $Im(\lambda (N^{n-2k}_{ext}))$ is inside the regular part
$\Sigma_{\hat h}^{reg} \subset \Sigma_{\hat h}$.

Let us denote by  $L^{62}_{cycl} \subset L^{62}$ the submanifold
(with  boundary) given by the formula $L^{62}_{cycl}= L^{62}
\cap U_{\Sigma^{reg}}$.
The mappings $\kappa_a$, $\rho_a$ are extendable
from $U_{\Sigma^{reg}}$ to $L^{62}_{cycl} \subset L^{62}$. Let us prove that these
mappings are extendable to  mappings $\kappa_a: L^{62} \to
K(\Z/2,1)$, $\rho_a: L^{62} \to K(\Z/4,1)$.

 The complement of thе submanifold $L^{62}_{cycl} \subset L^{62}$ is
denoted by $L^{62}_{\I_3} = L^{62} \setminus L^{62}_{cycl}$. The
submanifold $L^{62}_{\I_3}$ is a submanifold in the regular
$\varepsilon$--neighborhood of $h(X) \subset \R^n$. Obviously, the
structure group of the $\Z/2 \int \D_4$--framing of the normal
bundle of the manifold (with boundary) $L^{62}_{\I_3}$ is reduced to the subgroup
$\I_3 \subset \Z/2 \int \D_4$.

Let us consider  the mapping of  pairs $\mu_a \times \kappa_a:
(L^{62}_{cycl},
\partial L^{62}_{cycl}) \to (K(\Z/4,1) \times K(\Z/2,1),K(\Z/2,1)
\times K(\Z/2,1))$. Let us consider the natural projection $\pi_b:
\I_3 \to \I_b$. The extension of the mapping $\mu_a \times
\kappa_a$ to the required mapping $L^{62} \to K(\Z/4,1) \times
K(\Z/2,1)$ is given by the composition $L^{62}_{\I_3} \to
K(\I_3,1) \stackrel{\pi_{b,\ast}}{\longrightarrow} K(\I_b,1)
\subset K(\Z/4,1) \times K(\Z/2,1)$, where $\kappa_1 \in
K(\I_b;\Z/2)$ determines the inclusion $K(\I_b,1) \subset
K(\Z/2,1) \subset K(\Z/4,1)$.

Let us formulate the results in the following lemma.

\subsubsection*{Lemma 4}

--1. Let $n \ge 2^{13}-2$ and $k$, $n-4k=62$ satisfy the
conditions of Theorem 1 (in particular, an arbitrary element in
the group $Imm^{sf}(n-k,k)$ admits a retraction of the order $62$.
Then for arbitrarily small positive numbers $\varepsilon_1$,
$\varepsilon_2$, $\varepsilon_1 >> \varepsilon_2$ (the numbers
$\varepsilon_1$, $\varepsilon_2$ are the calibers of the regular
deformations in the construction of the $PL$--mapping $\hat h
:\hat X \to \R^n$ and of the immersion $g: N^{n-2k}
\looparrowright \R^n$ correspondingly) there exists the mapping
$m_a =(\kappa_a \times \mu_a): \Sigma_{h}^{reg} \to K(\Z/4,1)
\times K(\Z/2,1)$ under the following condition. The restriction
$m_a \vert_{\partial \Sigma_h^{reg}}$ (by $\partial
\Sigma_h^{reg}$ is denoted the part of the singular polyhedron
consisting of points on the diagonal) has the target $K(\Z/2,1)
\times K(\Z/2,1) \subset K(\Z/4,1) \times K(\Z/2,1)$ and is
determined by the cohomological classes $\kappa_{\hat X,1},
\kappa_{\hat X,2}$.

--2. The mappings $\kappa_a$, $\mu_a$ induces a mapping $(\mu_a
\times \kappa_a): L^{62} \to K(\Z/4,1) \times K(\Z/2,1)$ on the
self-intersection manifold of the immersion $g$.
\[  \]

Let us prove that the mapping $(\mu_a \times \kappa_a)$
constructed in  Lemma 4 determines a $\Z/2 \oplus
\Z/4$--structure for the $\D_4$--framed  immersion $g$.
We have to prove the equation (9).

Let us recall that the component $L^{62}_{int}$ of the
self-intersection manifold of the immersion $g$ is a $\Z/2 \int
\D_4$--framed manifold with  trivial Kervaire invariant: the
corresponding element in the group $Imm^{\Z/2 \int \D_4}(62,
n-62)$ is in the image of the transfer homomorphism. Therefore  it is sufficient to prove the equation
$$ <m_a^{\ast}(\rho \tau^{15} t^{31});[L^{62}]>=\Theta, $$
or, equivalently, the equation
$$ <(\hat \rho_a^{31} \hat \kappa_a^{31});[\hat L^{62}]>=\Theta, \eqno(10)$$
where  $\hat L \to L$  is the canonical cover over the
self-intersection manifold, $\hat L \subset N^{n-2k}_{ext}$
is the canonical inclusion.

By Herbert's theorem (see [A1] for the analogous construction)
we may calculate the right side of the equation by the formula
$$<\eta^{\ast}(w_2(\I_b))^{\frac{n-2k}{2}};[N^{n-2k}_{ext}/\sim]>. \eqno(11)$$
In this formula  by $N^{n-2k}_{ext}/\sim$ is denoted the quotient
of the boundary $\partial N^{n-2k}_{ext}=N^{n-2k-1}_Q$ that is
contracted onto the polyhedron $Q$ with the loss of the dimension.
Note that the mapping
 $m_a \vert_{N^{n-2k-1}_Q}$ is obtained by the composition of the mapping
$p_Q: N^{n-2k-1} \to Q$ with a loss of dimension  with the mapping $Q \to K(\I_b,1)$, the
last mapping is determined by the cohomology classes
$\kappa_{i,Q} \in H^1(Q;\Z/2)$, $i=1,2$. Therefore,
$m_{a\ast}([N^{n-2k}_{ext}/\sim]) \in H_{n-2k}(\I_b;\Z/2)$ is a
permanent cycle and the integration over the cycle
$[N^{n-2k}_{ext}/\sim]$ of the inverse image of the universal
cohomology class in (11) is well-defined.

It is convenient to consider the characteristic number  $\Theta_a$
as the value of a homomorphism $H_{n-2k}(X;\Z/2) \to \Z/2$ on the
cycle  $\lambda_{\ast}[N^{n-2k}_{ext}/\sim] \in H_{n-2k}(X;\Z/2)$.
This homomorphism is the result of the calculation of the
characteristic class $w_2(\I_b) \in H^2(K(\I_b,1);\Z/2)$ on the
prescribed cycle, i.e. on the image of the fundamental cycle
$[N^{n-2k}_{ext}/\sim]$ with respect to the mapping
$N^{n-2k}_{ext}/\sim \to \hat X \to K(\I_b,1)$. The cycle
$\lambda_{\ast}[N^{n-2k}_{ext}/\sim] \in H_{n-2k}(X;\Z/2)$ is the
modulo 2 reduction of an integral homology class. Therefore this cycle is given by a sum of fundamental
classes of the product of the two odd-dimensional projective
spaces, the sum of the dimensions of this spaces being equal to
$n-2k$.

Let us consider an arbitrary submanifold
 $S^{k_1}/i \times \RP^{k_2} \subset X$,
$k_1 + k_2 = \frac{n}{2}+31$, $k_1, k_2$ being odd. Let us consider
the cover $\RP^{k_1} \times \RP^{k_2} \to S^{k_1}/i \times
\RP^{k_2}$ and the composition $\RP^{k_1} \times \RP^{k_2} \subset
\hat X \stackrel{\hat h}{\looparrowright} \R^n$ after an
$\varepsilon_1$--small generic perturbation.  Let us denote this
mapping by $s_{k_1,k_2}$.

The self-intersection manifold of the generic mapping
$s_{k_1,k_2}: \RP^{k_1} \times \RP^{k_2} \to \R^n$ is a manifold
with boundary denoted by   $\Lambda^{62}_{k_1,k_2}$. The mapping
$$ \mu_a \times \kappa_a : (\Lambda^{62}_{k_1,k_2}, \partial N^{n-2k}_{k_1,k_2}) \to
(K(\Z/4,1)\times K(\Z/2,1), K(\Z/2,1) \times K(\Z/2,1))$$ is
well-defined. The $61$-dimensional homology fundamental class
$[\partial \Lambda]$ is integral, therefore the image of this
fundamental class $(\mu_a \times \kappa_a)_{\ast}([\partial
\Lambda^{62}_{k_1,k_2}]) \in H_{61}(K(\Z/4,1) \times
K(\Z/2,1);\Z/2)$ is trivial for a dimensional reason.

Therefore the homology class
$$(\mu_a \times \kappa_a)_{\ast}([\Lambda^{62}_{k_1,k_2}, \partial \Lambda^{62}_{k_1,k_2}])
\in $$ $$H_{62}(K(\Z/4,1)\times K(\Z/2,1),K(\Z/2,1) \times
K(\Z/2,1);\Z/2)$$ is well-defined. Let us consider the (permanent)
homology class
$$(\mu_a \times \kappa_a)^{!}_{\ast}([\bar \Lambda^{62}_{k_1,k_2}])
\in H_{62}(K(\Z/2,1)\times K(\Z/2,1);\Z/2), \eqno(12)$$ defined
from the relative class above by the transfer homomorphism.

To prove (10) it is sufficient to prove that the class (12)
coincides with the characteristic class
$$ p_{\ast, b} \circ \hat \eta_{\ast}([\hat \Lambda]) \in H_{62}(K(\I_b,1);\Z/2)$$
under the following isomorphism of the target group $\I_b = \Z/2
\oplus \Z/2$. By this isomorphism the prescribed generators in
$H^1(\Z/2 \oplus \Z/2;\Z/2)$ are identified with the cohomology
classes $\tau_{1}, \tau_2 \in H^1(K(\I_b,1);\Z/2)$  (compare with
Lemma 8 in [A1]). Theorem 2 is proved.

\section{Kervaire Invariant One Problem}

In this section we will prove the following theorem.

\subsubsection*{Main Theorem}

There exists an integer $l_0$ such that for an arbitrary integer $l \ge l_0$, $n=2^l-2$ the Kervaire invariant given by the formula (1) is trivial. 
\[  \]

\subsubsection*{Proof of Main Theorem}

Take the integer $k$ from the equation $n-4k=62$. Consider the diagram (5).
By the Retraction  Theorem [A2], Section 8 there exists an integer $l_0$ such that for an arbitrary integer $l \ge l_0$
an arbitrary
element $[(f,\Xi,\kappa)]$ in the 2-component of the cobordism group
$Imm^{sf}(\frac{3n+q}{4},\frac{n-q}{4})$ admits a retraction of order $62$. By Theorem 2
in the cobordism class $\delta[(f,\Xi,\kappa)]$ there exists a $\D_4$-framed immersion $(g,\Psi,\eta)$
with an $\I_4$-structure. 

Take the self-intersection manifold $L^{62}$ of $g$ and let $L_0^{10} \subset L^{62}$
be the submanifold dual to the cohomology class $\kappa_a^{28}\mu_a^{\ast}(\tau)^{12} \in H^{52}(L^{62};\Z/2)$.
By a straightforward calculation the restriction of the normal bundle of $L^{62}$ to the submanifold 
$L_0^{10} \subset L^{62}$ is trivial and the normal bundle of $L_0^{10}$ is the Whitney sum $12 \kappa_a \oplus 12 \mu_a$, where $\kappa_a$ is the line $\Z/2$-bundle, $\mu_a$ is the plane $\Z/4$-bundle with the characteristic
classes $\kappa_a$, $\mu_a^{ast}(\tau)$ described in the formula (8). By Lemma 6.1 (in the proof of this lemma we have to assume that the normal bundle of the manifold $L^{10}_0$ is as above) and by Lemma 7.1 [A2] the characteristic class (8) is trivial.
The Main Theorem is proved.

\[  \]
\[  \]

 Moscow Region, Troitsk, 142190, IZMIRAN.

 pmakhmet@mi.ras.ru
\[  \]

\end{document}